\documentclass[leqno,11pt,twoside]{article}
\usepackage{pliska}

\newcommand{\be}{\begin{equation}}
\newcommand{\ee}{\end{equation}}
\newcommand{\beq}{\begin{eqnarray}}
\newcommand{\eeq}{\end{eqnarray}}
\newcommand{\nbeq}{\begin{eqnarray*}}
\newcommand{\neeq}{\end{eqnarray*}}

\begin{document}

\head
{    
} {Critical Controlled Branching Processes and Their Relatives 
} { George P. Yanev  
} { Controlled Branching Processes  
} { George P. Yanev 
} {22
} {2015
}

\begin{abstract}
This survey aims at collecting and presenting results for one-type, discrete time branching processes with random control functions. In particular, the subclass of critical migration processes with different regimes of immigration and emigration is reviewed in detail. Critical controlled branching processes with continuous state space are also discussed.
\end{abstract}

\kwams{ 60J80 
}{ critical branching processes, immigration, emigration, control, extinction 
}

\bigskip

\normalsize \noindent 
\begin{tabular}[t]{l}
George P. Yanev\\
University of Texas - Pan American\\
yanevgp@utpa.edu
\end{tabular}

\section{Introduction}

The independence of individuals' reproduction is a fundamental assumption in branching processes. Since the 1960s, a number of authors have been studying models allowing different forms of population size dependence. Sevastyanov and Zubkov (1974) proposed a class of branching processes in which the number of reproductive individuals in one generation  decreases or increases depending on the size of the previous generation through a set of control functions. The individual reproduction law (offspring distribution) is not affected by the control and remains independent of the population size. These processes are known as controlled or $\varphi$-branching processes (CBP). N. Yanev (1975) (no relation to the author) extended the class of CBP by  introducing random control functions.
The so called $\varphi$-processes with random $\varphi$ can be defined as follows.

\vspace{0.3cm}{\bf Definition} The process  $\{Z_n, \ n\!=\!0,1,\ldots\}$ is called controlled branching process (CBP) if
\be \label{def1}
Z_{n+1}=\sum_{i\in I}\sum_{j=1}^{\varphi_{i,n}(Z_n)} \xi_{j,n}(i), \quad n\ge 0;\quad Z_0=z_0>0,
\ee
where $I$ is an (finite or infinite) index set and  for $i\in I$
\begin{description}
\item{(i)}
$\xi_i=\{\xi_{j,n}(i),\  j=1,2,\ldots;\ n=0,1,\ldots\}$
are i.i.d., non-negative, integer-valued r.v.'s, (independent for different $i$'s). Denote $\xi:=\xi_1$.
\item{(ii)}
$\varphi_i=\{\varphi_{i,n}(k),\   k=0,1, \ldots; \ n=0,1,\ldots\}$ are
  non-negative, integer-valued r.v.'s, independent from $\xi_i$, (independent
 for different $i$'s), and such that
$P(\varphi_{i,n}(k)=j)=p_k(j)$, for $j=0,1,\ldots$ Denote $\varphi(k):=\varphi_1(k)$.
\end{description}
The recurrence (\ref{def1}) describes a very large class of stochastic processes including, for  instance, all Markov chains with discrete time.
Among the particular cases of CBPs is the classical Galton-Watson process (GWP) as well as  popular discrete time branching processes such as: (i) processes with immigration:
$I=\{1,2\}$, $\varphi_{1,n}(k)=~k$, and  $\varphi_{2,n}(k)\equiv~1$; (iii) processes with state-dependent immigration:
$I=\{1,2\}$, $\varphi_{1,n}(k)=k$, and  $\varphi_{2,n}(k)=\max\{1-k,0\}$; and (iii) processes with random migration (to be discussed in Sections~3--5):
$I=\{1,2\}$, $\varphi_{1,n}(k)=\max\{\min\{k,k+\beta_n\},0\}$, and  $\varphi_{2,n}(k)=\max\{\beta_n,0\}$, where for $p+q+r=1$ we have $P(\beta_n=-1)=p$, $P(\beta_n=0)=q$, and $P(\beta_n=1)=r$.
In all these subclasses of CBPs, the controlled functions satisfy the condition
\[
\lim_{n\to \infty}\sum_{i\in I} \varphi_{i,n}(k)=\infty \quad \mbox{a.s.},\qquad k\ge 1,
\]
which can be identified as a general property of CBPs.

In Section~2 we present a classification of CBPs into subcritical, critical, and supercritical based on their mean growth rate. Two sets of conditions for extinction and non-extinction are presented  in relation to this classification. Finally, limit theorems for the critical CBPs are given. The next three sections, are devoted to critical processes with different regimes of migration.
In Section~3, processes with migration, stopped and non-stopped at zero, are defined and limit theorems in the critical case are discussed. Processes with time non-homogeneous migration are treated in Section~4. In Section~5, a more general type of migration is considered, utilizing a regenerative construction. CBPs with continuous state space are discussed in Section~6. Finally, the paper ends with some concluding remarks and a list of  references.

\section{General Class of Controlled Branching Processes}
In this section we shall discuss a classification of CBPs, which is similar to that of the classical Galton-Warson processes. Then we will focus our attention on the critical case.

As it will become clear below, the asymptotic behavior of the CBPs depends crucially on the so-called mean growth rate. Following Bruss (1984), we define the mean growth rate per individual in a population with $k$ mothers by
\[
\tau_k:=k^{-1}E[Z_{n+1}\ |\ Z_n=k]= k^{-1}E[\varphi(k)]E[\xi].
\]
In the particular case of GWP, we have $\tau_k= E[\xi]$, i.e., the mean growth rate equals the offspring mean and remains constant for any $k$.

\subsection{Extinction and Classification of CBPs}

Extinction, along with growth and composition of the population, is a principal subject of interest in the theory of branching processes.
If the control functions satisfy $\varphi_n(0)\equiv 0$ a.s., then $\{Z_n\}$ is a Markov chain with absorption state 0. Furthermore, it can be proven (see \cite{Y75}) that if
$P(\xi=0)>0$ or $P(\varphi(k)=0)>0$ for $k=1,2,\ldots$,
then the classical extinction-explosion duality
\[
P(Z_n\to 0)+P(Z_n\to \infty)=1
\]
holds. The following key theorem for the extinction probability of CBPs is proven in Gonzales et al. (2002).

\vspace{0.3cm}{\bf Theorem 2.1} (\cite{GMP02})
\begin{description}
\item{(i)} If $\limsup_{k\to \infty}\tau_k< 1$, then $P(Z_n\to 0\ |\ Z_0=N)=1$ for $N\ge 1$.
\item{(ii)}
If
$
\liminf_{k\to \infty}\tau_k>1,
$
 then there exists $N_0$ such that for $N\ge N_0$ we have
 $P(Z_n\to 0\ |\ Z_0=N)<1$.
\end{description}
Referring to Theorem 2.1, Gonzalez et al. (2005) classify the CBPs as follows.

\vspace{0.3cm}{\bf Definition} The class of CBPs can be partitioned into three subclasses:
\begin{description}
\item{(i)} subcritical if $\limsup_{k\to \infty}\tau_k<1$;
\item{(ii)} critical if
$\liminf_{k\to \infty}\tau_k\le 1\le \limsup_{k\to \infty}\tau_k$;
\item{(iii)} supercritical if $\liminf_{k\to \infty}\tau_k>1$.
\end{description}
Unlike the supercritical GWP, if the number of ancestors in the supercritical CBP is not sufficiently large, then the extinction probability might be one. This resembles the situation with the two--sex branching processes (e.g., \cite{H93}). On the other hand, the critical CBP does not always have extinction probability one, as it is seen in (\ref{tau2}) below.

Gonzalez at al. (2005) study in detail the extinction probability of the critical CBP considering different rates of convergence of $\tau_k$ to one. Their findings rely on analysis of the stochastic difference equation
\be \label{SDE}
Z_{n+1}=Z_n+h(Z_n)+\delta_{n+1} \quad \mbox{a.s.}, \quad\qquad n=0,1,\ldots,
\ee
where $h(k)=E[\varphi(k)]E[\xi]-k$ and $\delta_{n+1}=Z_{n+1}-E[Z_{n+1}|Z_n]$ (martingale difference). Denote
\[
2\sigma_{s}(k):=E[|\delta_{n+1}|^s|Z_n=k]\qquad s>0.
\]
It is proven in \cite{GMP05} that if
\be \label{tau1}
\lim_{k\to \infty}\tau_k=1\qquad \mbox{and}\quad \tau(k)\ge 1,
\ee
then for $N\ge 1$
\be \label{tau2}
P(Z_n\to 0\ |\ Z_0=N)
    \left\{
\begin{array}{ll}
=1  & \mbox{if} \quad \limsup_{k\to \infty} (\tau_k-1)k^2(\sigma_2(k))^{-1}<1 ;\\
\\
<1 & \mbox{if} \quad \liminf_{k\to \infty} (\tau_k-1)k^2(\sigma_2(k))^{-1}>1.
   \end{array}
   \right.
\ee
A different set of conditions for extinction and non-extinction of $\{Z_n\}$ is obtained by  N. Yanev (1975) using a random walk construction.
It is proven in \cite{Y75} that if the control functions have a linear growth a.s., that is
\[
\varphi_n(k)=\alpha_nk\ (1+o(1))\quad \mbox{a.s.} \quad k\to \infty,
\]
where $\{\alpha_n\}$ are i.i.d. and independent of the reproduction, then for $N\ge 1$
\[
P(Z_n\to 0\ |\ Z_0=N)
    \left\{
\begin{array}{ll}
=1  & \mbox{if} \quad E[\log (\alpha_1 E\xi)]<0;\\
\\
<1 & \mbox{if} \quad  E[\log (\alpha_1 E\xi)]>0.
   \end{array}
   \right.
\]
Bruss (1980) shows that the independence of reproduction assumption for $\{\alpha_n\}$ can be removed.

\subsection{Limit Theorems for Critical CBPs}
Assuming (\ref{tau1}), let us turn to the critical CBPs.
It follows from (\ref{tau2}) that, depending on the rate
 of convergence of $\tau_k$ to one, the extinction probability is either one or less than one.
We will consider these two cases separately. Utilizing a gamma--limit theorem (see \cite{K92}) for the stochastic difference equation (\ref{SDE}),
Gonzalez et al. (2005) prove the following two theorems.

\vspace{0.3cm}{\bf Case A.} The extinction is almost sure, i.e.,
 $P(Z_n\to 0\ |\ Z_0=N)=1$.

\vspace{0.3cm}{\bf  Theorem 2.2} (\cite{GMP05}) 
Assume
\begin{description}
\item{(i)} $\tau_k=1+ck^{-1}, \quad c>0,\ $ $\ k=1,2,\ldots$;
\item{(ii)} $\sigma_2(k)=2ak+O(1), \quad a>0,\ $ as $\ k\to \infty$;
\item{(iii)} $\sup_{k\ge 1}\left(g_k^{1/k}\right)^{'''}(1)<\infty$, where $g_k(s):=E\left[ s^{\varphi(k)}\right]$, $0\le s\le 1$. 
\end{description}
If $c\le a$,  then
\be \label{exp_lim}
\lim_{n\to \infty} P\left(\frac{Z_n}{an}\le x|Z_n>0\right)=1-e^{-x}.
\ee
Note that the limiting distribution in (\ref{exp_lim}) is exponential as in
the critical GWP. However, one difference is that $P(Z_n>0)\sim c_1n^{-(1-c/a)}$,
 $c_1>0$ whereas the survival probability in the GWP has a
decay rate $(bn)^{-1}$.

\vspace{0.3cm}{\bf Case B.} Positive non-extinction probability, i.e., $P(Z_n\to 0\ |\ Z_0=N)<1$.

\vspace{0.3cm}{\bf  Theorem 2.3} (\cite{GMP05}) 
Assume as $k\to \infty$
\begin{description}
\item{(i)} $\tau_k=1+ck^{-(1-\alpha)}+o\left(k^{-(1-\alpha)}\right), \quad c>0, \ 0<\alpha<1$;
\item{(ii)} $\sigma_2(k)=2ak^{1+\alpha}+o\left(k^{1+\alpha}\right), \quad a>0$;
\item{(iii)} $\sigma_{2+s}(k)=O\left(\left(\sigma_2(k)\right)^{1+s/2}\right)\quad \mbox{for some}\ s>0$.
\end{description}
If $c>a$,  then
\be \label{gamma_lim}
\lim_{n\to \infty}P\left(\frac{Z_n^{1-\alpha}}{(1-\alpha)^2an}\le x\ |\ Z_n>0\right) = \frac{1}{\Gamma(\gamma)}\int_0^x t^{\gamma-1}e^{- t}\, dt,
\ee
where
$\gamma=(c-a\alpha)/(a(1-\alpha))$ and $\Gamma(x)$ is the Gamma function. Obviously,
if $\alpha=0$ and $c=a$, then (\ref{gamma_lim}) coincides with (\ref{exp_lim}).

\section{Branching Processes with Migration}

In the context of queueing theory, stochastic models with migration
 were discussed in \cite{K79}.
A model of branching process with emigration-immigration (migration) was introduced in Nagaev and Han (1980). The readers are referred to the survey \cite{VZ93}
 and the monograph \cite{R95} for a detailed account of results for processes with a variety of immigration and emigration regimes.

In Sections~3--5 we shall review results for a class of CBPs, called branching processes with migration, introduced in
N. Yanev and Mitov (1980) when a detailed study of these branching models began. These processes were already mentioned in Section~1 as a special class of CBPs.
The particular choice of control functions $\varphi(k)$ allows for a detailed analysis which in turn leads to interesting new findings.
On the other hand, branching processes with migration are sufficiently general to include as subclasses previously studied models with different regimes of immigration and emigration.

\vspace{0.3cm}{\bf Definition} The process $\{Y_n,\  n=0,1,\ldots\}$ is called a branching process with migration if $Y_0>0$ and for $n=0,1,\ldots$
\be \label{BPM}
Y_{n+1}=\left\{
     \begin{array}{ll}
\sum_{k=1}^{Y_n} \xi_{k,n}+M^+_n  & \mbox{if} \quad Y_n>0;\\
\\
M^0_n  & \mbox{if} \quad Y_n=0,
   \end{array}
   \right.
\ee
where for $p+q+r=1$
\be \label{M+}
M^+_n=\left\{
       \begin{array}{lll}
    -\xi_{1,n}   & \mbox{probab.} \quad p,  &  \quad \mbox{(emigration)}  \\
    \ \  0                & \mbox{probab.} \quad q,  & \quad \mbox{(no migration)}\\
    \ \  \eta_n & \mbox{probab.} \quad r, &  \quad \mbox{(immigration)}
        \end{array}
   \right.
\ee
is the migration outside zero and
\[
M^0_n=\left\{
       \begin{array}{lll}
  \ \  0                & \mbox{probab.} \quad 1-r, & \quad \mbox{(no migration)} \\
\ \ \eta_n  & \mbox{probab.} \quad r, &  \quad \mbox{(immigration at 0)}, \\
       \end{array}
   \right.
\]
is the migration at zero. The number of immigrants $\{\eta_n, \ n=1,2,\ldots\}$ are i.i.d. non-negative, integer-valued, and independent from the offspring variables.
The process $\{Y_n\}$ can be interpreted as follows. Three scenarios are possible: (i) the offspring of one individual is removed (emigration) with probability $p$; (ii) there is no migration with probability $q$; or (iii) $\eta_n$ individuals join the population (immigration) with probability $r$. The state zero is a reflecting barrier for $\{Y_n\}$. The emigration here can be regarded as "reversed" (negative) immigration since the branching process is modified to allow both positive and negative increments. See \cite{P14} and the references within for other approaches to emigration.

\subsection{Branching Migration Processes with Reflection Barrier at Zero}
 In Sections~3--5 we assume (unless stated otherwise) that $\{Y_n\}$ is critical with finite offspring variance and finite immigration mean, i.e.,
\be \label{main_ass}
E[\xi]=1, \qquad  2b:=Var [\xi]<\infty, \qquad \mbox{and}\qquad d:=E[\eta_n]< \infty.
\ee
The long-term behavior of the critical $\{Y_n\}$ depends crucially on the parameter
\be \label{theta_def}
\theta:=\frac{EM^+_n}{(Var[\xi])/2}=\frac{rE[\eta_n]-pE[\xi]}{(Var[\xi])/2}=\frac{rd-p}{b},
\ee
i.e., the ratio of the mean migration outside zero over half of the offspring variance.
Depending on the values of $\theta$, the aperiodic and irreducible  Markov chain $\{Y_n\}$ can be classified as
\[
\{Y_n\}=
     \left\{
        \begin{array}{ll}
     \mbox{non-recurrent} &  \quad \ \ \, \theta > 1 \\
     \mbox{null-recurrent} & 0\le \theta \le 1\\
     \mbox{positive-recurrent} & \quad \ \ \, \theta< 0.
         \end{array}
     \right.
\]

\noindent The following limiting results are obtained in \cite{YY96}.

\vspace{0.3cm}{\bf  Theorem 3.1} (\cite{YY96}) Assume (\ref{main_ass}).
\begin{description}
\item{(A)} If $\theta>0$ (dominating immigration) and $Var[\eta_n]<\infty$, then
\[
   \lim_{n\to \infty}P\left(\frac{Y_n}{bn}\le x\right) =
\frac{1}{\Gamma(\theta)}\int_0^x t^{\theta-1}e^{-t}\, dt.
\]
\item{(B)} If $\theta=0$ (zero average migration), then
\[
\lim_{n\to \infty}P\left(\frac{\log Y_n}{\log n} \le x\right) = x, \qquad 0<x<1.
\]
\end{description}

\begin{description}
\item{(C)} If $\theta<0$ (dominating emigration), then there is a limiting-stationary distribution, i.e.,
\[
\lim_{n\to \infty}P(Y_n=k) = v_k, \qquad \sum_{k=0}^\infty v_k=1
\]
and $V(s)=\sum_{k=0}^\infty v_ks^k$ is the unique p.g.f. solution of a functional equation.
\end{description}

{\bf Remarks} (i) Theorem 3.1(C) is similar in form to a theorem for processes with immigration due to Seneta (1968). Both theorems claim the existence of a limiting-stationary distribution in the critical case.
 One obvious difference between the two results is that the former applies to the migration process with
dominating emigration, i.e., $\theta<0$ whereas the latter is for a process with immigration. Another important difference is that Theorem 3.1(C) assumes
finite offspring variance whereas in Seneta's theorem the offspring variance is assumed infinite.

(ii) It is worth pointing out here a limit theorem due to Dyakonova (1997) for the "close to critical" process $\{Y_n\}$, i.e., assuming that the offspring mean $m:=E[\xi]\uparrow 1$. It is known that if $m<1$, then $\{Y_n\}$ has a limiting-stationary distribution. Let $V$ be the limiting random variable with this distribution. Then it is proven in \cite{D97} that
\[
\lim_{m\uparrow 1} P\left( \frac{\log V}{\log \frac{1}{1-m}}\le x\right)=x, \qquad x\in(0,1).
\]

\subsection{Branching Migration Processes with Absorbing State Zero}
\vspace{0.3cm}In this subsection we shall consider the branching process with no migration when it hits zero,
i.e., $M^0_n =0$ a.s. in (\ref{BPM}) and hence zero is an absorbing state.

\vspace{0.3cm}{\bf Definition} Let $Y_0^0>0$ and for $n=1,2,\ldots$
\[
Y_n^0=Y_n I_{\{Y_n>0\}} \quad \mbox{a.s.},
\]
where $I_A$ denotes the indicator of the event $A$. Then $\{Y_n^0\}$ is called a migration process with absorption at zero.

It is shown in \cite{YY95} and \cite{YY04}, under some additional finite moment conditions, that the probability of
the process surviving to time $n$ satisfies as $n\to \infty$
\[
P(Y^0_n>0)\sim
     \left\{
        \begin{array}{ll}
      c_\theta>0 &  \quad \ \ \, \theta > 1 \\
     c_\theta\  (\log n)^{-1}& \quad \ \ \, \theta =1\\
     c_\theta \ n^{-(1-|\theta|)} & \  \! 0\le \theta<1,\\
     c_\theta\  n^{-(1+|\theta|)} & \quad \ \ \, \theta<0.
         \end{array}
     \right.
\]
Referring to these results, one can adopt the following classification for the critical process $\{Y^0_n\}$: (i) critical-supercritical for $\theta>1$; critical-critical for $0\le \theta \le 1$; and (iii) critical-subcritical for $\theta< 0$.

The next two limit results were proven in \cite{YVM86} (for $\theta>0$) and
\cite{YY95} (for $\theta\le 0$).

\newpage\vspace{0.3cm}{\bf Theorem 3.2} (\cite{YVM86}, \cite{YY95}) 
Assume (\ref{main_ass}).
\begin{description}
\item{(A)}
If $\theta>1$ (strongly dominating immigration), then
\[
\lim_{n\to \infty} P\left( \frac{Y^0_n}{bn}\le x|Y^0_n>0\right)=
\frac{1}{\Gamma(\theta)}\int_0^x t^{\theta-1}e^{-t}\, dt.
\]
\item{(B)}
 Assume $\theta\le 1$ and some additional moment conditions when $\theta<0$ (see \cite{YY95}).
Then
\[
\lim_{n\to \infty} P\left( \frac{Y^0_n}{bn}\le x|Y^0_n>0\right)=1-e^{-x}.
\]
\end{description}

{\bf Remarks} (i) If the rate of migration is not too high, i.e., $\theta\le 1$, then the long-term behavior of $\{Y^0_n\}$ over
the non-extinction trajectories is the same as in the critical GWP. The observation after (\ref{exp_lim}) applies here too.

(ii) One extension of Theorem~3.2(A), when the distribution of the initial number of ancestors $Y^0_0$ belongs to the domain of attraction
of a stable law with parameter in $(0,1]$, is given in \cite{YY97}.

\section{Time Non-Homogeneous Migration}

N. Yanev and Mitov (1985) study branching processes with time non-homogeneous migration defined as follows.

\vspace{0.3cm}{\bf Definition} The process $\{\tilde{Y}_n: n=0,1,\ldots\}$ is called a branching process with non-homogeneous migration if $\tilde{Y}_0>0$ and for $n=1,2,\ldots$
\be \label{BPMn}
\tilde{Y}_{n+1}=\left\{
     \begin{array}{ll}
\sum_{k=1}^{\tilde{Y}_n} \xi_{k,n}+\tilde{M}^+_n  & \mbox{if} \quad \tilde{Y}_n>0;\\
\\
\tilde{M}^0_n  & \mbox{if} \quad \tilde{Y}_n=0,
   \end{array}
   \right.
\ee
where the migration is given for $p_n+q_n+r_n=1$ by
\be \label{M+n}d
\tilde{M}^+_n=\left\{
       \begin{array}{lll}
    -\xi_{1,n}   & \mbox{probab.} \quad p_n,  &  \quad \mbox{(emigration)}  \\
    \ \  0                & \mbox{probab.} \quad q_n,  & \quad \mbox{(no migration)}\\
    \ \  \eta_n& \mbox{probab.} \quad r_n, &  \quad \mbox{(immigration)}
        \end{array}
   \right.
\ee
 and
\[
\tilde{M}^0_n=\left\{
       \begin{array}{lll}
  \ \  0                & \mbox{probab.} \quad 1-r_n, & \quad \mbox{(no migration)} \\
\ \ \eta_n  & \mbox{probab.} \quad r_n, &  \quad \mbox{(immigration at 0)}. \\
       \end{array}
   \right.
\]
Unlike (\ref{M+}), here the probabilities $p_n, q_n$ and $r_n$ controlling the migration are time-dependent. Thus, $\{\tilde{Y}_n\}$ is a non-homogeneous Markov chain.
In addition to ({\ref{main_ass}), suppose that the immigration variance is finite, i.e.,
\be \label{assump nonhom}
Var[\eta_n]<\infty.
\ee
In the rest of this section we also assume that the migration decreases to 0, i.e., $\lim_{n\to \infty}q_n=~1$.

\vspace{0.3cm}{\bf Case A.} Decreasing to Zero Migration and $ p_n=o(r_n)$.

\vspace{0.3cm}{\bf Theorem 4.1} (\cite{YM84}, see also \cite{YM85a}) Suppose (\ref{main_ass}) and (\ref{assump nonhom}). If as $n\to \infty$
\[
r_n\sim \frac{r}{\log n}, \qquad  p_n=o(r_n),
\]
 then
\[
\lim_{n\to \infty}P\left(\frac{\log Z_n}{\log n}\le x \right)=  e^{-rd(1-x)/b} \quad 0\le x\le 1.
\]

\vspace{0.3cm}{\bf Theorem 4.2} (\cite{YM85b}, see also \cite{YM85a}) Suppose (\ref{main_ass}) and (\ref{assump nonhom}). If as $n\to\infty$
\[
r_n\sim \frac{l_n}{\log n}\quad \mbox{and}\quad p_n=o(r_n),
\]
where  $l_n\sim o(\log n)\to\infty$, then for $x\ge 0$
\[
\lim_{n\to \infty}P\left(l_n\left(1-\frac{\log Z_n}{\log n}\right)\le x\right)
=1-e^{-dx/b}.
\]

\vspace{0.3cm}{\bf Case B.} Decreasing to Zero Migration and $ p_n=dr_n$.

\vspace{0.3cm}If both immigration and emigration decrease to zero at the same rate, then a key role for the limiting behavior of the process is played by
the series $\sum_{n=0}^\infty p_n$ and  $\sum_{n=0}^\infty r_n$. This observation is made precise in the next three theorems.


\vspace{0.3cm}{\bf Theorem 4.3} (\cite{YM85a}, \cite{DLY06})
Suppose
(\ref{main_ass}), (\ref{assump nonhom}), and $p_n=dr_n$. If one of the following two conditions holds as $n\to \infty$

(i) $p_n\sim l_n n^{-v}$ for $0<v<1$ where $l_n$ is a s.v.f. at $\infty$.

(ii) $p_n= O\left((\log n)^{-1}\right)$,

\noindent then
\[
\lim_{n\to \infty}P\left(\frac{ \log \tilde{Y}_n}{\log n}\le x|\tilde{Y}_n>0\right)= x,
\quad x\in (0,1).
\]
Theorem 4.3 is an analog of Foster's result for
processes with immigration at zero only. Unlike Foster's model, $\{\tilde{Y}_n\}$ is a non-homogeneous Markov chain.

\vspace{0.3cm}{\bf Theorem 4.4} (\cite{YM85a}) Suppose
(\ref{main_ass}), (\ref{assump nonhom}), and $p_n=dr_n$.  Assume as $n\to \infty$
\begin{description}
\item{(i)} $p_n \sim l_nn^{-1}$, where $l_n$ is a s.v.f. at $\infty$;
\item{(ii)} $\lim_{n\to \infty} \frac{p_nn\log n}{\sum_{k=1}^n p_k}=C$ for
    $0\le C\le \infty$.
\end{description}
Then,
\be \label{P1}
\lim_{n\to \infty}P\left(\frac{\log \tilde{Y}_n}{\log n}\le x|\tilde{Y}_n>0\right)= \frac{C}{1+C}x=:G_1(x)
\quad 0<x<1
\ee
and
\be \label{P2}
\lim_{n\to \infty}P\left(\frac{ \tilde{Y}_n}{bn}\le x|\tilde{Y}_n>0\right)=
\frac{C}{1+C} + \frac{1}{1+C}\left(1-e^{-x}\right)=:G_2(x)
\quad x>0.
\ee
It is worth pointing out that, since $\lim_{x\to 1}G_1(x)=\lim_{x\to 0}G_2(x)$, the limiting distributions in  (\ref{P1}) and (\ref{P2}) represent the two different types of non-degenerate trajectories of $\{\tilde{Y}_n\}$:
\begin{description}
\item{(A)} $\tilde{Y}_n\sim n^{\eta_1}$, where $\eta_1\in U(0,1)$ with probab. $\frac{C}{1+C}$;
\item{(B)} $\tilde{Y}_n\sim \eta_2n$, where $\eta_2\in Exp(b)$ with probab. $\frac{1}{1+C}$.
\end{description}
We will have a similar situation with the processes considered in Section~6.

\vspace{0.3cm}{\bf Theorem 4.5} (\cite{YM85a}) Suppose
(\ref{main_ass}), (\ref{assump nonhom}), and $p_n=dr_n$.  If $\sum_{k=1}^\infty p_k<\infty$, then
\[
\lim_{n\to \infty}P\left(\frac{ \tilde{Y}_n}{bn}\le x|\tilde{Y}_n>0\right)=
1-e^{-x},
\quad x\ge 0.
\]
Theorem 4.5 is an analog of the classical Kolmogorov--Yaglom result for GWPs. It turns out that the convergence of $\sum_{n=0}^\infty p_n$ and $\sum_{n=0}^\infty r_n$ ensure that the migration disappears without a trace so fast that the process with non-homogeneous migration has the same asymptotic behavior as
the standard GWP.

\section{Regenerative Branching  Processes with Migration}

Quoting \cite{W14}, "A regenerative process is a stochastic process with the property that after some (usually) random time, it
starts over in the sense that, from the random time on, the process is stochastically equivalent to what it
was at the beginning". Regenerative processes can be intuitively seen as comprising of i.i.d. cycles. For classical regenerative processes, cycles and cycle lengths are i.i.d.

 Consider a random vector $(W,R)$ with non-negative and independent coordinates. The sequence of its i.i.d. copies $(W_j,R_j)$ for $j=1,2,\ldots$ defines an alternating renewal process (e.g., \cite{MO14}).
 The random variables $W$ and $R$ can be interpreted as the "working" and "repairing" time periods, respectively, of an operating system. Denote $S_0=0$ and for $n=1,2,\ldots$
\[
 S_n:=\sum_{j=1}^n (R_j+W_j) \qquad \mbox{and}\qquad
N(t):=\max\{n\ge 0:\ S_n\le t\}.
\]
Define
\[
\sigma(t):=t-S_{N(t)}-R_{N(t)+1}, \qquad t\ge 0.
\]
The random variable $\sigma(t)$ takes on positive or negative values depending on whether
at $t$ the system is working or repairing, respectively.
Let associate with each $W_j$, $j\ge 1$ a cycle given by the process
$\{Z_j(t), \ 0\le t\le W_j\}$ such that
\[
Z_j(0)=0, \quad Z_j(t)>0\ \mbox{for}\ 0<t<W_j, \quad Z_j(W_j)=0.
\]

\vspace{0.3cm}{\bf Definition} An alternating regenerative process (ARP) is defined by
\[
Z(t):= \left\{
     \begin{array}{ll}
          Z_{N(t)+1}(\sigma(t)) & \mbox{when}\  \sigma(t)\ge 0\quad \mbox{(the system is working)}\\
            0   & \mbox{when} \  \sigma(t)<0\quad \mbox{(the system is repairing)}.
      \end{array}
       \right.
\]

\vspace{0.3cm}{\bf Example} Recall the process with migration $\{Y_t,\ t=0,1,\ldots\}$ defined by (\ref{BPM}). It is an ARP with $P(R=k)=\left[P(M^0_t=0)\right]^{k-1}\left[1-P(M^0_t=0)\right]$,
a geometrically  distributed repairing time.
Consider the sequence $\{Y^0_{t,j}, \ j=1,2,\ldots\}$ of corresponding migration processes with absorption at 0 and let $W_j=tI_{\{Y^0_{t,j}>0\}}$ for $j\ge 0$.
Thus, $\{Y_t\}$ is an ARP with cycle process $\{Y^0_{t,j}\}$. It regenerates whenever it visits state zero.

The  migration process in the example above can be generalized as follows.

\vspace{0.3cm}{\bf Definition} Define
 a  regenerative branching process with migration by $X_0=~0$ and for $t=1,2,\ldots$
\[
X_t= \left\{
     \begin{array}{ll}
          Y^0_{N(t)+1,\sigma(t)} & \mbox{when}\  \sigma(t)\ge 0\\
            0   & \mbox{when} \  \sigma(t)<0,
      \end{array}
       \right.
\]
where $\{Y^0_{j,t}, \ j=1,2,\ldots\}$ are migration processes with absorbing state zero.

Note that, unlike the process with migration $\{Y_t\}$, in the generalized regenerative
branching  process with migration $\{X_t\}$ the repairing time periods $R_j$ are not necessary geometrically distributed.

\vspace{0.3cm}{\bf Possible Scenario.} The queueing systems are good examples for discrete time regenerative processes. Consider a single-server queue with Poisson arrivals. The service periods are
composed of a busy part (not--empty queue) $W_j$ and an idle part (empty queue) $R_j$. The customers arriving during the service time of a customer are her "offspring". The "immigrants" (probably from another customer pool) will be served in the end of the entire "generation". Alternatively,
some "emigrants" may give up and leave the queue.

\vspace{0.3cm}Let $\{X_t\}$ be critical and $0<\theta<1/2$, where $\theta$ is from (\ref{theta_def}).
Assume that either  $E[R]$ is finite or $P(R>t)\sim L(t)t^{-\alpha}$ for $\alpha\in (1/2,1]$, where $L(t)$ is a s.v.f.
Under some additional moment assumptions for the reproduction and migration,
the following limiting results are obtained in  G. Yanev, Mitov, and N. Yanev (2006). The proofs make use of theorems from Mitov and N. Yanev (2001) for regenerative processes.

(i) If "the working time dominates over the repairing time", i.e.,
\[
0\le c:=\lim_{t\to\infty}\frac{P(R>t)}{P(W>t)}<\infty,
\]
then for $x\ge 0$
\[
\lim_{t\to \infty}P\left(\frac{X_t}{bt}\le x\right)= \frac{c}{c+1}+
\frac{1}{c+1}\frac{1}{B(\theta, 1-\theta)}\int_0^1 y^{\theta-1}(1-y)^{-\theta}\left(1-e^{-x/y}\right)\, dy,
\]
where $B(x,y)$ is the Beta function. The expected value of the limiting random variable is $\theta/(c+1)$.

(ii) If  "the repairing time dominates over the working time", i.e.,
\[
\lim_{t\to\infty}\frac{P(R>t)}{P(W>t)}=\infty,
\]
then for $x\ge 0$
\[
\lim_{t\to \infty}P\left(\frac{X_t}{bt}\le x\ |\ X_t>0\right)=
\frac{1}{B(\theta, \alpha)}\int_0^1 y^{\theta-1}
(1-y)^{\alpha-1}\left(1-e^{-x/y}\right)\, dy.
\]
Note that the distribution of the limiting random variable is a mixture of beta and exponential distributions and has a mean of $\theta/(\theta+\alpha)$.

\section{Controlled Branching Processes with Continuous State Space}

A branching process with continuous state space models situations when it is difficult to
count the number of individuals in the population, but a related non--negative variable (e.g., volume or weight) associated with the "individuals" is measured instead.

\vspace{0.2cm}Let us make the following assumptions.
\begin{description}
\item{(i)} For fixed $n$, let $U_n:=\{U_{i,n}, \ i\ge 1\}$ be a sequence of i.i.d., non-negative random variables and the double array $U:=\{U_n, \ n\ge 1\}$ consists of independent
sequences $U_n$, n=1,2,\ldots
\item{(ii)} Each of the stochastic processes $N_n:=\{N_n(t), \ t\in T\}$, n=1,2,\ldots has state space $Z^+$, the set of non-negative integers. They are independent processes with stationary and independent increments (s.i.i.) and $N_n(0)=0$ a.s. Here $T$ is either $[0,\infty)$ or $Z^+$.
\item{(iii)} The sequence $V:=\{V_n, \ n\ge 1\}$ consists of independent and  non-negative random variables.
\item{(iv)} The processes $N:=\{ N_n,\ n\ge 1\}$, $U$, and $V$ are independent.
\item{(v)} The random variable $X_0$ is non-negative and independent from all processes introduced in (i)-(iv).
\item{(vi)} The components $N_n$ and $U_n$ for $n\ge 1$ of the processes $N$ and $U$, respectively, are identically distributed.
\end{description}
 The following class of branching processes is introduced by Adke and Gadag (1995) and studied by Rahimov (2007) and Rahimov and Al--Sabah~(2008).

\vspace{0.3cm}{\bf Definition}  A controlled branching process with continuous state space
is defined by the recursive relation
\be \label{def8}
X_{n+1}=\sum_{i=1}^{N_{n+1}(X_n)} U_{i,n+1} + V_{n+1}, \quad n=0,1,\ldots,\quad X_0=0.
\ee
Notice that if $X_0$, $U$, and $V$ are integer-valued, then $\{X_n\}$ is a CBP. If, in addition, we choose in (\ref{def1}) the index set to be $I=\{1,2\}$ and the control functions to be $\varphi_{1,n}(k)=N_n(k)$, and $\varphi_{2,n}(k)\equiv 1$, we obtain (\ref{def8}).

It is proven in \cite{AG95} that $\tilde{Z}_n:=N_n(X_{n-1})$ for $n=1,2,\ldots$ is a GWP with time-depended immigration given by
\[
\tilde{Z}_{n+1}=\sum_{i=1}^{\tilde{Z}_n} \xi_{i, n+1}+\eta_{n+1},\quad n=0,1,\ldots,\quad \tilde{Z}_0=0,
\]
where $\xi_{i,n+1}\stackrel{d}{=}N_{n+1}(U_n)$ and
$\eta_{n+1}\stackrel{d}{=}N_{n+1}(V_n)$. Exploring this duality, Rahimov (2007) transferred results from GWPs with immigration to $\{X_n\}$. Below we present one theorem from \cite{R07} for the critical $\{X_n\}$. Denote
 $2\tilde{b}:=Var[\xi_{i,n}]$, $\beta_n:=E[\eta_n(\eta_n-1)]$, and $\gamma_n:=E[V_n]$. For simplicity, some of the assumptions of the next theorem are given in terms of moments of $\xi_{i,n}$ and $\eta_n$, which
 can be expressed as functions of the moments of $N$, $U$, and $V$ (see \cite{R07}).

\vspace{0.3cm}{\bf Theorem 6.1} (\cite{R07}) 
 Suppose $E[N_1(1)]E[U_1]=1$ and
$\tilde{b}< \infty$. Assume
\begin{description}
\item{(i)} $\beta_n=o\left(\gamma_n\log n\right)\to 0$ as $n\to \infty$;
\item{(ii)} $\lim_{n\to \infty}\gamma_n\log n=0$ and $\lim_{n\to \infty} \frac{\gamma_nn\log n}{\sum_{k=1}^n \gamma_k}=C$ for
    $0\le C\le \infty$.
\end{description}
Then,
\[
\lim_{n\to \infty}P\left(\frac{\log X_n}{\log n}\le x|X_n>0\right)= \frac{C}{1+C}x,
\quad 0<x<1
\]
and
\[
\lim_{n\to \infty}P\left(\frac{ X_n}{
\tilde{b}n}\le x|X_n>0\right)=
\frac{C}{1+C} + \frac{1}{1+C}\left(1-e^{-x}\right)
\quad x>0.
\]
The similarities between Theorem 4.2 and Theorem 6.1 are striking. The phenomenon of having different limiting distributions under different normalization in GWPs with decreasing  time-dependent immigration was observed by Badalbaev and Rahimov (1978) (see also \cite{R95}, p.109 and p.122).

\section{Concluding Remarks}

This survey is by no means exhaustive. Not included here are some of classes CBPs such as:
branching processes with barriers (see Zubkov (1972), Bruss (1978), Schuh (1976), Sevastyanov (1995)), CBPs with random environments,
and the more recently introduced alternating branching processes (see Mayster (2005)).

Controlled branching processes are part of Sevastyanov's legacy. Over
time, particular subclasses were introduced and studied in details. We paid special attention to the processes with migration,
which have been a subject of systematical research investigations by the Bulgarian school in branching processes under the direction of its founder Professor Nikolay Yanev a.k.a. the Captain.
Closed relations were established between CBP and other classes, e.g.,
two-sex processes and population size-dependent processes.
There is no doubt, that CBPs have great potential as modeling tools.
In my opinion, they deserve more attention from the branching processes'
community.


\begin{thebibliography}{9999}
\bibitem{AG95} \textsc{S. R. Adke and V. G. Gadag.}  A new class of branching processes. Branching Processes. In: Ed. Heyde, C. C., Proceedings of the First World Congress,
(Lecture Notes in Statistics 99), 90--105,  New York, Springer,  1995.
\bibitem{BR78} \textsc{I. S. Badalbaev and I. Rahimov.} Limit theorems for critical Galton-Watson processes with immigration decreasing intensity.\textit{ Izv. Akad. Nauk UzSSR, Ser. Fiz.- Mat. Nauk}, \textbf{89}(1978) 2:9--14,(Russian).
 \bibitem{BR78b} \textsc{I. S. Badalbaev and I. Rahimov.} Critical branching processes with decreasing intensity.\textit{Theory Prob. Appl.}, \textbf{23}(1978) 2:259--268.
\bibitem{B78} \textsc{F. T. Bruss.} Branching processes with random absorbing processes. \textit{J. Appl. Prob.}, \textbf{15}(1978), 54--64.
\bibitem{B80} \textsc{F. T. Bruss.} A counterpart of the Borel-Cantelli lemma. \textit{ J. Appl. Prob.}, \textbf{17}(1980), 1094--1101.
\bibitem{B84} \textsc{F. T. Bruss.} A note on extinction criteria for bisexual Galton-Watson processes.\textit{ J. Appl. Probab.}, \textbf{21}(1984), 915--919.
\bibitem{DLY06} \textsc{M. Drmota, G. Louchard, N. M.  Yanev.} Analysis of a recurrence related to critical nonhomogeneous branching processes. \textit{Stoch. Anal. Appl.}, \textbf{24}(2006), 37--59.
\bibitem{D97} \textsc{E. E. Dyakonova.} Close-to-critical branching processes with migration. \textit{Theory Probab. Appl.}, \textbf{41}(1997), 151--156.
\bibitem{GMP02} \textsc{M. Gonzalez, M. Molina, and I. del Puerto.} On the class of controlled branching processes with random control functions. \textit{J. Appl. Prob.}, \textbf{39}(2002), 804--815.
\bibitem{GMP05} \textsc{M. Gonzalez, M. Molina, and I. del Puerto.} Asymptotic behaviour of critical controlled branching processes with random control functions.
\textit{ J. Appl. Prob.}, \textbf{42}(2005), 463--477.
\bibitem{H93} \textsc{D. M. Hull.} How Many Mating Units Are Needed to Have a Positive Probability of Survival? \textit{Mathematics Magazin},
\textbf{66}(1993), 28-33.
\bibitem{K79} \textsc{F. P. Kelly} Reversibility and Stochastic Networks. New York, Wiley, 1979.
\bibitem{K86} \textsc{G. Kersting.} On recurrence and transience of growth models. \textit{J. Appl. Prob.}, \textbf{23}(1986), 614--625.
\bibitem{K92}  \textsc{G. Kersting.} Asymptotic $\Gamma$-distribution for stochstic difference equations. \textit{Stoc. Proc. Appl.}, \textbf{40}(1992), 15--28.
\bibitem{M01} \textsc{P. Mayster.} Alternating branching processes. \textit{J. Appl. Probab.}, \textbf{42}(2005), 1095--1108.
\bibitem{MO14} \textsc{K. Mitov and E. Omey.} Renewal Processes. New York, Springer, 2014.
\bibitem{MY01} \textsc{K. Mitov and N. Yanev.}  Limit theorems for alternating renewal processes in the infinite mean case. \textit{Adv. in Appl. Probab.}, \textbf{33}(2001), 896–911.
\bibitem{NH80}  \textsc{S. V. Nagaev and L. V. Han.} Limit theorems for critical Galton-Watson branching process with migration. \textit{Theory Prob. Appl.}, \textbf{25}(1980), 3:523--534.
\bibitem{P14} \textsc{A. G. Pakes.} Immigration-Emigration Processes. \textit{Wiley StasRef: Statistic Reference Online}, pp.5, 2014.
\bibitem{R95} \textsc{I. Rahimov.} Random Sums and Branching Stochastic Processes. (LNS, 96), New York, Springer, 1995.
\bibitem{R07} \textsc{I. Rahimov.} Asymptotic behavior of a controlled branching process with continuous state space. \textit{Stoc. Anal. Appl.}, \textbf{25}(2007), 337--352.
\bibitem{R08} \textsc{I. Rahimov and W. S. Al-Sabah.} Limiting behavior of a generalized branching process with immigration. \textit{Stats. Probab. Lett.}, \textbf{78}(2008), 225--230.
\bibitem{VZ93} \textsc{V. A. Vatutin and A. M. Zubkov.} Branching Processes II. \textit{J. Soviet. Math.}, \textbf{67}(1993), 6.
\bibitem{S95} \textsc{B. A. Sevastyanov.} Extinction probabilities of branching processes bounded from below. \textit{Theory Probab. Appl.}, \textbf{40}(1995), 495-–502.
\bibitem{SZ74} \textsc{B. A. Sevastyanov and A. M. Zubkov.} Controlled branching processes. \textit{Theory Prob. Appl.}, \textbf{19}(1974), 15--25.
\bibitem{S76} \textsc{H. -J. Schuh.} A condition for the extinction of a branching process with an absorbing lower barrier. \textit{J. Math. Biology}, \textbf{3}(1976), 271--287.
\bibitem{W14} \textsc{R. W. Wolff.} Regenerative Processes. \textit{Wiley StasRef: Statistic Reference Online}, pp.7, 2014.
\bibitem{YMY06} \textsc{G. Yanev, K. Mitov and N. Yanev.} Critical branching regenerative process with random migration. \textit{J. Appl. Statist. Sci.}, \textbf{12}(2006), 41--54.
\bibitem{YY95} \textsc{G. Yanev and N. Yanev.}  Critical branching process with random migration. In: Ed. C.C. Heyde, Branching Processes,  Lecture Notes in Statistics 99,  New York, Springer, 36--46, 1995.
\bibitem{YY96} \textsc{G. Yanev and N. Yanev.} Branching processes with two types emigration and state-dependent immigration. In: Eds. Heyde, C. C., Prohorov, Yu. V., Pyke, R., Rachev, S. T. Athens Conference on Applied Probability and Times Series,  Vol. 1, Applied Probability, Lecture Notes in Statistics 114, New York, Springer, 216--228, 1996.
\bibitem{YY97} \textsc{G. Yanev and N. Yanev.} Limit theorems for branching processes with random migration stopped at zero. In: Eds. K.B. Athreya and P. Jagers,  Classical and Modern Branching Processes, The IMA Volumes in Math. and Its Appl., Vol. 84, New York, Springer, 323--336, 1997.
\bibitem{YY04} \textsc{G. Yanev and N. Yanev.} A critical branching process with stationary-limiting distribution. \textit{Stoch. Anal. Appl.}, \textbf{22}(2004), 721--738.
\bibitem{Y75} \textsc{N. Yanev.} Conditions for degeneracy of $\varphi$- branching processes with random $\varphi$. \textit{Theory Prob. Appl.}, \textbf{20}(1975), 433--440.
\bibitem{YM80}  \textsc{N. Yanev and K. Mitov.} Controlled branching processes: the case of random migration. \textit{C.R. Acad. Bulg. Sci.}, \textbf{33}(1980), 473--475.
\bibitem{YM84}\textsc{N. Yanev and K. Mitov.} Controlled branching processes with non-homogeneous migration.\textit{ Pliska Stud. Math. Bulgar.}, \textbf{7}(1984), 90--96, (Russian).
\bibitem{YM85a} \textsc{N. Yanev and K. Mitov.} Critical branching processes with nonhomogeneous migration. \textit{Ann. Probab.}, \textbf{13}(1985), 3:923--933.
\bibitem{YM85b} \textsc{N. Yanev and K. Mitov.} A critical branching process with decreasing migration.\textit{ Serdica (Bulg. Math. Publ.)}, \textbf{11}(1985), 240--244.
\bibitem{YVM86} \textsc{N. Yanev, V. Vatutin, and K. Mitov.} \textit{Math. and Math. Edu.}, 1986, 511--517 (Russian).
\bibitem{Z72}  \textsc{A. M. Zubkov.} A degeneracy condition for a bounded branching process. \textit{Mat. Zametki}, \textbf{8}(1970), 9--18, (Russian).
\bibitem{Z74} \textsc{A. M. Zubkov.} Analogies between Galton-Watson processes and $\varphi$-branching processes. \textit{Theory Prob. Appl.}, \textbf{19}(1974), 319--339.
\end{thebibliography}
\end{document}